\documentclass[reqno]{amsart}

\usepackage{enumerate}
\usepackage{ifpdf}
\usepackage{amsmath}
\usepackage{amsfonts}
\usepackage{amssymb}
\usepackage{amsthm}
\usepackage{hyperref}
\usepackage{setspace}
\usepackage{amsrefs}

\newcommand{\ignore}[1]{}

\renewcommand{\Re}{\operatorname{Re}}
\renewcommand{\Im}{\operatorname{Im}}

\newcommand{\abs}[1]{\left\lvert {#1} \right\rvert}
\newcommand{\norm}[1]{\left\lVert {#1} \right\rVert}

\newcommand{\C}{{\mathbb{C}}}
\newcommand{\R}{{\mathbb{R}}}

\newcommand{\bB}{{\mathbb{B}}}

\newcommand{\sA}{{\mathcal{A}}}

\newcommand{\sH}{{\mathcal{H}}}

\newcommand{\sO}{{\mathcal{O}}}

\newcommand{\sY}{{\mathcal{Y}}}

\newcommand{\rank}{\operatorname{rank}}

\newtheorem{thm}{Theorem}[section]

\newtheorem{prop}[thm]{Proposition}

\newtheorem{lemma}[thm]{Lemma}
\newtheorem{claim}[thm]{Claim}

\theoremstyle{definition}

\theoremstyle{remark}

\theoremstyle{remark}
\newtheorem{example}[thm]{Example}

\author{Ji\v{r}\'i Lebl}
\address{Department of Mathematics, University of Illinois
at Urbana-Champaign, 
Urbana, IL 61801, USA}
\email{jlebl@math.uiuc.edu}
\date{July 7, 2008}

\ifpdf
\fi

\title
{Levi-flat hypersurfaces with real analytic boundary}

\begin{document}


\begin{abstract}
Let $X$ be a Stein manifold of dimension at least 3.  Given a compact
codimension 2 real analytic submanifold $M$ of $X$, that is the boundary of a
compact Levi-flat hypersurface $H$, we study the regularity of $H$.  Suppose
that the CR singularities of $M$ are an $\sO(X)$-convex set.  For example,
suppose $M$ has only finitely many CR singularities, which is a generic
condition.  Then $H$ must in fact be a real analytic submanifold.  If $M$ is
real algebraic, it follows that $H$ is real algebraic and in fact extends past
$M$, even near CR singularities.  To prove these results we provide two
variations on a theorem of Malgrange, that a smooth submanifold contained in a
real analytic subvariety of the same dimension is itself real analytic.  We
prove a similar theorem for submanifolds with boundary, and another one for
subanalytic sets.
\end{abstract}

\maketitle


\section{Introduction} \label{section:intro}

Let $X$ be a Stein manifold of dimension $N$.
Suppose that $M$ is a compact real submanifold of $X$ of codimension 2.
We ask the following question:
\emph{Does there exist a compact Levi-flat hypersurface $H$ with boundary
$M$, and if so what is the regularity of $H$?}  In this paper we address the
regularity and uniqueness of $H$ in the case $M$ is real analytic.
In particular, if the CR singularities of $M$ are an $\sO(X)$-convex set,
$H$ must in fact be a real analytic submanifold.
Generically a real submanifold $M$ of codimension 2
has only finitely many CR singularities, which of course
is an $\sO(X)$-convex set.
It is also possible to interpret these results as a regularity and uniqueness
statement about a certain nonlinear partial differential equation.
In \cite{Lebl:ext}
we saw that the real analytic case is rigid locally, and hence the
present results are a natural continuation.

Malgrange \cite{Malgrange} proved that if
a smooth submanifold is contained in a real analytic subvariety of the same
dimension then the submanifold is real analytic.  In \S \ref{section:varmalg}
we prove two variations of this theorem which are of independent interest.
One for subanalytic sets, and one for submanifolds with
boundary.

We will suppose that $N \geq 3$.  When $X$ is two dimensional,
there is no CR geometric 
information on $M$, and the methods presented here do not apply.
In particular, when $N=2$, any real analytic foliation
on a real analytic
submanifold of codimension 2 locally extends to a Levi-flat hypersurface.
On the
other hand a codimension 2 real submanifold
is not in general even locally the boundary of a Levi-flat hypersurface.

For a two-dimensional $X$, similar questions, both locally and globally,
have been studied for many years.
For example by Bishop \cite{Bishop:diffman},
Moser and Webster \cite{MW:normal}, Bedford \cite{bedford:boundaries}, 
Huang and Krantz \cite{huangkrantz}, or
Bedford and Gaveau \cite{BG:envhol}, among others.
Our question can be interpreted as a form of a complex Plateau problem,
see Bedford \cite{bedford:boundaries} for example.
The existence question for $N\geq 3$ has been considered by
Dolbeault, Tomassini and Zaitsev \cite{DTZ:CRmath}.


First we fix some terminology.
In the sequel, \emph{submanifold} $M \subset X$ will always mean an embedded
submanifold, i.e.\@ the topology of $M$ is the subspace topology induced by
$X$, though $M$ need not be closed in $X$.
\emph{Hypersurface} will mean a submanifold of codimension 1, and
\emph{real analytic subvariety} of an open set $U$,
a set closed in $U$ and locally defined by
the vanishing of a family of real analytic functions.
We will say $V \subset \R^N$ is a \emph{local real analytic subvariety}
if there exists an open set $U$ such that $V$ is a real analytic
subvariety of $U$.
A submanifold will be called \emph{real algebraic} if it is contained
in a real algebraic subvariety (a subvariety defined by
the vanishing of real polynomials) of the same dimension.

We give the following definitions assuming $X = \C^N$ for clarity.  The
definitions are local in character and carry over to an arbitrary Stein
manifold in a straightforward way.

Let $M \subset \C^N$ be a real submanifold.
Let $J$ be the complex structure on $\C^N$, and let
$T_p^cM = J(T_pM) \cap T_pM$.
$M$ is said to be \emph{CR} near $p$ if there exists a neighbourhood $U$
of $p$ such that $\dim T_q^cM$ is constant as $q$ varies in $U$.
The points near which $M$ is CR will be denoted by $M_{CR}$, and the
complement of $M_{CR}$ we will call the CR singular
points.
The set of CR singular points is a real analytic subvariety of $M$.

A set $H \subset \R^n$ is
a \emph{$C^k$ hypersurface with boundary}, if
there is a subset $\partial H \subset H$,
such that
$\partial H \subset H$,
$H \setminus \partial H$ is a $C^k$ hypersurface
(submanifold of codimension 1), and for each point $p \in
\partial H$,
there exists a neighbourhood $p\in U \subset \R^n$,
a $C^k$ diffeomorphism $\varphi \colon U \to \R^n$,
such that
$\varphi (H \cap U) = \{ x \in \R^n \mid x_{n-1} \geq 0, ~ x_{n} = 0\}$,
and such that
$\varphi (\partial H \cap U) = \{ x \in \R^n \mid x_{n-1} = 0,~ x_{n} =
0\}$.  Hence, $\partial H$ is a $C^k$ submanifold of codimension 2
in $\R^n$.
We will call $H^o := H \setminus \partial H$ the \emph{interior} of $H$.
A $C^k$ ($k \geq 2$) hypersurface $H \subset \C^N \cong \R^{2N}$ is said to be \emph{Levi-flat}
if the bundle $T^cH$ is involutive.
If $H$ is a hypersurface with boundary, then 
we will say it is Levi-flat when $H^o$ is Levi-flat.

Our first result is the following.  For a Stein manifold $X$,
let $\sO(X)$ denote the holomorphic
functions on $X$.  For a compact set $S$,
let $\hat{S}$ denote the $\sO(X)$-convex hull of $S$.

\begin{thm} \label{regthm}
Let $X$ be a Stein manifold of dimension $N \geq 3$.
Let $M \subset X$
be a compact real analytic submanifold of codimension 2.
Suppose there exists a compact connected $C^\infty$ Levi-flat
hypersurface $H$ with boundary,
such that $\partial H = M$.

Let $S$ be the set of CR singularities of $M$.
If $S = \hat{S}$, then $H^o = H \setminus M$
is real analytic.
Further,
$H$ is the unique compact connected Levi-flat $C^\infty$ hypersurface with
boundary $M$.
\end{thm}

In particular if $S$ is a finite set, then $S = \hat{S}$ and
$H^o$ is real analytic and unique.
Since $M$ is an oriented codimension 2 compact submanifold,
it is a standard result that 
the condition of only isolated CR singularities (hence finitely
many) is the generic situation.  We need only look at the Gauss
map of the manifold and apply Thom's transversality theorem.
See Lai \cite{Lai:immerse} for example.

We will see that $H$ actually extends as a real analytic submanifold
past all the CR points, though at CR singularities the picture is not
clear.

Let us interpret this result as a regularity and uniqueness statement for
a certain partial differential equation.
Suppose that $\rho$ is a function with $d\rho \not= 0$.
The hypersurface $\{ \rho = 0 \}$
is Levi-flat if and only if the Levi form vanishes
\cites{BER:book, Boggess:CR, DAngelo:CR}.  That is equivalent to the complex
Hessian vanishing on all holomorphic vectors tangent to the hypersurface.
Hence $\{ \rho = 0 \}$ 
is Levi-flat if and only if the complex bordered Hessian
is of rank two on the hypersurface.  In other words,
$\{\rho = 0 \}$ is Levi-flat if and only if
\begin{equation} \label{flateqnice}
\rank
\left[
\begin{matrix}
\rho & \rho_z \\
\rho_{\bar{z}} & \rho_{z\bar{z}}
\end{matrix}
\right]
= 2
\ \ \ \text{ for all points on $\{\rho = 0\}$. }
\end{equation}
As \eqref{flateqnice} implies the determinant of the complex bordered Hessian
is zero,
the equation is related to equations of the complex Monge-Amp{\`e}re
type.

We can also think of the hypersurface as a graph of a certain function.
This leads us to the following somewhat more complicated differential equation.
Let $\Omega \subset \C^{N-1} \times \R$ be a domain and let us
call the coordinates $(z,s) \in \C^{N-1} \times \R$, and suppose
$u \colon \Omega \to \R$ is a function that satisfies
\begin{equation} \label{flateq}
\begin{split}
& u_{z_j} (-i+u_s) u_{s\bar{z}_j}
+
u_{\bar{z}_j} (i+u_s) u_{sz}
+
\\
& \ \ \ \ \ \ \ 
u_{z_j} u_{\bar{z}_j} u_{ss}
-
(1+u_s^2)u_{{z_j}\bar{z}_j} = 0
\ \ \ \text{ for all $j=1, \ldots, N-1$},
\\
&
u_{z_j} u_{\bar{z}_k} u_{z_k\bar{z}_j}
+
u_{\bar{z}_j} u_{z_k} u_{z_j\bar{z}_k}
-
\\
& \ \ \ \ \ \ \ 
u_{z_j} u_{\bar{z}_j} u_{z_k\bar{z_k}}
-
u_{z_k} u_{\bar{z}_k} u_{z_j\bar{z}_j}
= 0
\ \ \ \text{ for all $j,k=1, \ldots, N-1$}.
\end{split}
\end{equation}
We can consider the graph of this function to lie in
$\C^N$ by looking at the set defined by $\Im w = u(z,\bar{z},\Re w)$.
The graph is then a Levi-flat hypersurface, as
$\rho(z,\bar{z},w,\bar{w}) = \Im w - u(z,\bar{z},\Re w)$
satisfies \eqref{flateqnice}.

Now consider the following boundary value problem.
Suppose that $\Omega$ has real analytic boundary.
Let 
$g \colon \partial \Omega \to \R$ be a real analytic function,
and suppose there exists a solution $u \in C^\infty (\overline{\Omega})$,
that satisfies \eqref{flateq} and $g = u|_{\partial \Omega}$.
Further, impose the generic condition
that at most at finitely many
points is $\partial \Omega$ tangent to a line $\{ s = s_0 \}$, for some
$s_0 \in \R$.
Then Theorem \ref{regthm} tells us that $u$ is real analytic
on $\Omega$.  Further, $u$ is the unique solution in
$C^\infty(\overline{\Omega})$.

Let us now look at Theorem \ref{regthm} when $M$
is a real algebraic submanifold.
In this case,
we will be able to extend the hypersurface $H$
past $M$ even near CR singularities, see Lemma \ref{algreglocal}.
Using this lemma we extend Theorem \ref{regthm} as follows.

\begin{thm} \label{algreg}
Let $M \subset \C^N$, $N \geq 3$,
be a compact real algebraic submanifold of codimension 2.
Suppose there exists a compact connected $C^\infty$ Levi-flat
hypersurface $H$ with boundary,
such that $\partial H = M$.

If the set of CR singularities of $M$ is
polynomially convex, then
there exists a real algebraic
Levi-flat hypersurface $\sH$ (without boundary) such that $H \subset \sH$.
\end{thm}

Theorems \ref{regthm} and \ref{algreg} do not hold as stated
when $N=2$.  For the uniqueness, see Example
\ref{example:uniq}.  For more discussion and examples on the
regularity in two dimensions see the Appendix in
\cite{bedford:boundaries}.

More background details about CR geometry and some of the methods employed here
can be found in the books
\cites{BER:book, Boggess:CR, DAngelo:CR}.

In \S \ref{section:varmalg} we discuss two variations on
a theorem of Malgrange that will be needed in the proofs.
In \S \ref{section:regthm}
we prove Theorem \ref{regthm}.
In \S \ref{section:algreg}
we discuss the real algebraic case and prove
Theorem \ref{algreg}.
In \S \ref{section:hypersing} we consider hypersurfaces with isolated
singularities.
In \S \ref{section:examples} we look at some examples.

The author would like to acknowledge Peter Ebenfelt for many useful discussions
while preparing these results.  Also the author would like to acknowledge
John D'Angelo, Edward Bierstone, and Dmitri Zaitsev
for their useful comments on the results and the manuscript.
The author would like to acknowledge Alexander Tumanov for suggesting
to consider singular hypersurfaces as is done in section
\S \ref{section:hypersing}.

\section{Two variations on the theorem of Malgrange} \label{section:varmalg}

Malgrange proved
that a smooth submanifold contained in a real analytic subvariety of the same
dimension is itself real analytic.
See \cite{BER:book} Theorem 5.5.32 or \cite{Malgrange} Chapter VI,
Proposition 3.11.
We will prove two variations on this theorem, one for subanalytic sets,
and one for submanifolds with boundary.
As these results are of independent interest,
we will prove these theorems for arbitrary codimension even though we will
only apply them in the hypersurface case.
The first theorem we wish to prove is a subanalytic version of Malgrange's
theorem.
The class of {\it semianalytic} sets is the smallest class containing
sets of the form $\{ f > 0 \}$ for real analytic functions $f$ that is
closed under finite union, intersection and complement.
A set is {\it subanalytic} if it is locally a projection of a
relatively compact semianalytic set.
See Bierstone and Milman \cite{BM:semisub} for more information
about subanalytic and semianalytic sets.

\begin{thm} \label{thm:subansmooth}
Suppose that $Y \subset \R^N$
is a $C^\infty$ submanifold,
and $S \subset \R^N$ a subanalytic set
of same dimension as $Y$, such that $Y \subset S$.  Then $Y$ is real analytic.
\end{thm} 

Next we want to prove a version for submanifolds with boundary.

\begin{thm} \label{thmmalvar}
Suppose that
$Y \subset \R^N$
is a $C^\infty$ submanifold with boundary,
and $V \subset \R^N$ a local real analytic subvariety
of same dimension as $Y$, such that $Y \subset V$.  Then there
exists a real analytic submanifold $\sY$ of same dimension
as $Y$ such that $Y \subset \sY$.
\end{thm}

In particular,
such $Y$ extends uniquely past its boundary as a real analytic
submanifold.
Combining the ideas of the proofs we can see that we could 
replace ``local real analytic subvariety'' in Theorem \ref{thmmalvar}
with ``subanalytic set.''

As we will be concerned with convergence we will require the
following classical theorem.  A proof can be found for example in
\cite{BER:book} (Theorem 5.5.30).

\begin{thm} \label{convofseriesthm}
Suppose $T(x)$ is a formal power series in $x \in \R^N$.  Suppose 
$T(tv)$ is a convergent power series in $t \in \R$ for
all $v \in \R^N$.  Then $T$ is convergent.
\end{thm}

In particular, one corollary of this theorem is that if $f(x)$ is a
$C^\infty$ function, and there exists a neighborhood $0 \in U \subset \R^N$
such that $t \mapsto f(t x)$ is real analytic whenever $t x \in U$,
then $f$ is real analytic.

\begin{proof}[Proof of Theorem \ref{thm:subansmooth}]
Let $Y$ and $S$ be as in the statement of Theorem \ref{thm:subansmooth}.
First suppose that $Y$ is a hypersurface.
Without loss of generality,
let $(x,y)\in \R^m \times \R$ be our coordinates and suppose
that $Y$ is defined by $y = f(x)$ for a $C^\infty$ function $f$,
with $f(0) = 0$.
We need only to prove that $Y$ is real analytic near $0$.

We restrict to a line, and
look at the graph of the function $t \mapsto f(t x)$.
The graph is an intersection of $Y$ and a $2$ dimensional plane $P_x$.
Suppose that $S \cap P_x$ was of dimension 2.  Let
$\tilde{S} = \overline{S \setminus P_x}$ which is again a subanalytic set
containing $Y$ and of the same dimension as $Y$.
A subanalytic set is a locally finite union of real analytic submanifolds.
Further, there is a sequence $\{ x_j \}$ converging to $x$ such
that $\dim S \cap P_{x_j} = 1$.  Hence, $\tilde{S} \cap P_x$ must be
of dimension 1.

Hence, the graph of $t \mapsto f(t x)$ is contained in a subanalytic
set of dimension 1.
By a theorem of {\L}ojasiewicz (see \cite{BM:semisub} Theorem 6.1),
a one dimensional subanalytic set is semianalytic, and hence contained
in a real analytic subvariety of the same dimension.
We can now apply the standard version of the theorem of Malgrange
to see that $t \mapsto f(t x)$
must be real analytic.  By the discussion after
Theorem \ref{convofseriesthm},
we know that $f$, and therefore $Y$ itself must be real analytic near $0$.

Now suppose that $Y$ is of higher codimension.  That is,
let $(x,y)\in \R^m \times \R^n$ be our coordinates and suppose
that $Y$ is defined by $y = f(x)$ for a $C^\infty$ mapping $f$,
with $f(0) = 0$.  Write $f = (f_1,\ldots,f_n)$.  Pick
$1 \leq j \leq n$, and let $\tilde{Y}$ be the submanifold
defined by
$y_j = f_j(x)$ in $\R^m \times \R$.
Take $\pi$ to be the projection $(x,y) \mapsto
(x,y_j)$.  As $S$ has dimension $m$, $\pi(S)$ is a subanalytic set of
dimension at most
$m$.  Because $\pi(S)$ contains $\tilde{Y}$,
it must be of dimension $m$ exactly.
Applying the result for hypersurfaces, we see that $f_j$ is real
analytic near $0$ for all $j$.  Therefore, $Y$ is real analytic near $0$.
\end{proof}

To prove Theorem \ref{thmmalvar},
it is enough to show that near every point $p \in Y$, $Y$ is
contained in a real analytic submanifold of same dimension.
The techniques used in this proof will be similar to those
in \cite{BER:book}.
If $p$ is not on the boundary, the statement follows from the 
standard statement of the theorem of Malgrange.
What we will prove, therefore, is the case
when $p \in \partial Y$.  We can assume that near $p$, the boundary
of $Y$ is real analytic, since we can always find a smaller submanifold
$\tilde{Y} \subset Y$ of equal dimension and with real analytic boundary
such that $p \in
\tilde{Y}$.  To find $\tilde{Y}$ we intersect $Y$ and $V$ with a
real analytic nonsingular hypersurface $B$ transversal to $Y$ such that
$B \cap Y$ is a smooth submanifold through $p$ and is contained in $B \cap V$,
which is of same dimension as $B \cap Y$.  Take $B \cap Y$ be the new
boundary which is real analytic by the standard statment of the theorem.
If $\tilde{Y}$ is real analytic at $p$, then it uniquely extends as a real
analytic submanifold past $p$.  $Y$ must therefore be contained in this
extension by analytic continuation, since we know that $Y$ is real analytic
in the interior.

We will assume that $p$ is the origin.
We first prove the codimension 1 case.  Assume that $Y$ is a hypersurface
with boundary.
After straightening the boundary
of $Y$, applying the
Weierstrass preparation theorem to $V$ at 0, and rescaling, we see that it suffices to prove the following lemma.
We will denote by $\bB_N \subset \R^N$
the unit ball centered at 0.

\begin{lemma} \label{genmalkeylemma}
Put
$\Omega = \bB_N \cap \{ x \in \R^N \mid x_N \geq 0 \}$.
Let
\begin{equation}
P(x,y) = y^k + c_{k-1}(x) y^{k-1} + \cdots + c_0(x),
\end{equation}
where $c_j$ are real analytic functions on $\bB_N$,
and $c_j(0) = 0$ for $j=0,\ldots,k-1$.
Suppose $f \colon \Omega \to \R$ is a $C^\infty$ function such that
\begin{equation}
f(0) = 0, ~\text{ and }~
P(x,f(x)) = 0 ~ \text{ for $x \in \Omega$}.
\end{equation}
Then $f$ is real analytic near 0.
\end{lemma}

What must be shown is that the formal Taylor series of $f$ converges
for $x$ near 0, and converges to $f$ whenever such $x$ is in $\Omega$.
We will want to apply Theorem \ref{convofseriesthm}, and so
we will want to reduce the proof of
Lemma \ref{genmalkeylemma} to the one dimensional case.
Here we will use the following statement of the Puiseux theorem.  A proof 
can be found for example in \cite{BM:arcanal} or \cite{BER:book}.

\begin{thm}[Puiseux]
Let $P(x,y) = y^k + c_{k-1}(x) y^{k-1} + \cdots + c_0(x)$,
where $c_j \in \C\{x\}$ and $c_j(0) = 0$ for $j=0,\ldots,k-1$.
Suppose that $P$ is irreducible, then there exists $a \in \C\{x\}$ such that
\begin{equation} \label{puiseuxeq}
P(t^k,y) = \prod_{j=1}^k \big(y-a(\omega^j t)\big) ,
\end{equation}
where $\omega = e^{i 2\pi / k}$.
\end{thm}

In one dimension, Lemma \ref{genmalkeylemma} reduces to the following
statement.

\begin{lemma} \label{onedimcase}
Let $P(x,y) = y^k + c_{k-1}(x) y^{k-1} + \cdots + c_0(x)$,
where $c_j$ are real analytic functions on $(-1,1)$ such that $c_j(0) = 0$.
Suppose $f \colon [0,1) \to \R$ is a $C^\infty$ function such that
$f(0) = 0$, and
$P(x,f(x)) = 0$ for $x \in [0,1)$.  Then $f$ is real analytic.
\end{lemma}

\begin{proof}
By the theorem of Malgrange, we know that $f$ is real analytic on
$(0,1)$.  Hence the graph of $f$ must lie in a single branch of the
zero locus of $P$, and we can assume that $P$ is irreducible.
If we plug in $f(t^k)$ for $y$ into \eqref{puiseuxeq}, we see that there exists
a $j$ such that $f(t^k) = a(\omega^j t)$, for $t \geq 0$.  Hence the
Taylor series of $t \mapsto f(t^k)$ at 0 is equal to the 
power series expansion of $t \mapsto a(\omega^j t)$.
Therefore the Taylor series
of $f$ also converges, and converges to $f$.
\end{proof}

\begin{proof}[Proof of Lemma \ref{genmalkeylemma}]
Let $T_f$ be the formal Taylor series
of $f$ at the origin.  Pick an arbitrary $v \in \R^N$, and
without loss of generality suppose that $v_N \geq 0$ and
$\norm{v} = 1$.
The function $g(t) = f(tv)$ is a $C^\infty$ function for $t \in [0,\epsilon)$
for some $\epsilon > 0$.  By Lemma \ref{onedimcase}, the Taylor series of
$g$ at 0 converges, and applying 
Theorem \ref{convofseriesthm} we see that $T_f$ converges as well.
Suppose it
converges for $\norm{x} < \delta$.
The $\delta$ therefore no longer depends on $x$.
Lemma \ref{onedimcase} tells us that the Taylor series of $g$
converges to $f(tv)$ for $t \in [0,\delta)$, and as it converges
for all $v$ as above, $T_f$ converges to $f$
for $\norm{x} < \delta$ and $x_N \geq 0$ as desired.
\end{proof}

To finish the proof of Theorem \ref{thmmalvar}, we only need to reduce to
the hypersurface case.  Let $Y$ and $V$ be as in the statement of the
theorem.  After a linear change of coordinates we only need
to apply the partial
generalization of Weierstrass Preparation Theorem to $V$ (see for
example \cite{BER:book} Theorem 5.3.9).  That is, we note that locally
near 0 there exists a subvariety $\tilde{V}$ of same
dimension as $V$, such that $V \cap U \subset \tilde{V} \cap U$ for
some small neighbourhood $U$ of 0, and such that $\tilde{V}$ is defined
in the coordinates $(x,y) \in \R^m \times \R^n$ by
\begin{equation}
y_k^{d_k} + \sum_{j=0}^{d_k-1} a_{jk}(x)y_k^j
\ \ \ \text{ for $k = n,\ldots,N$, }
\end{equation}
where $a_{jk}$ are real analytic functions vanishing at $0$.
Therefore, as we can assume that $Y$ is a graph over the $x$ coordinates,
we see that the theorem follows by applying Lemma
\ref{genmalkeylemma}.

\section{Real analytic regularity} \label{section:regthm}

The following more general, statement implies the
regularity part of Theorem \ref{regthm}.  Let $\sO(X)$ denote the holomorphic
functions on $X$.

\begin{thm} \label{regthmtech}
Let $M \subset X$
be a compact real analytic submanifold of codimension 2.
Suppose there exists a compact connected $C^\infty$ Levi-flat
hypersurface $H$ with boundary,
such that $\partial H = M$.

Let $S$ be the set of CR singular points of $M$, and $\hat{S}$
the $\sO(X)$-convex hull of $S$.  Then $H \setminus (M \cup \hat{S})$
is a real analytic submanifold.
\end{thm}

If a hypersurface $H$ is Levi-flat, then it is locally foliated by complex
hypersurfaces.  This foliation is called the \emph{Levi foliation}.
The smallest germ (in terms of dimension) of a
CR submanifold $W$ of $M$ through $p$ such that the $T_q^c W = T_q^c M$
is called the \emph{local CR orbit} at $p$, and is guaranteed to exist
by a theorem of Nagano
in case $M$ is real analytic, and by a theorem of Sussmann
in case $M$ is smooth.
The following two results have been in some form known previously.
Proofs using similar notation
as the present paper can be found in \cite{Lebl:ext}.

\begin{lemma} \label{lemma:nomin}
Suppose $H \subset \C^N$, $N \geq 2$,
is a $C^2$ Levi-flat hypersurface with boundary,
$M \subset \C^N$ is a $C^\infty$
submanifold of codimension 2, and $M = \partial H$.
If $M$ is CR near $p \in M$, then the local CR orbit of $M$
at $p$ is of positive
codimension in $M$ (i.e.\@ $M_{CR}$ is nowhere minimal).
\end{lemma}

\begin{lemma} \label{lemma:uniqueH}
Suppose $H \subset \C^N$, $N \geq 3$, is a $C^2$ Levi-flat hypersurface
with boundary,
$M \subset \C^N$ is a real analytic
submanifold of codimension 2, and $M = \partial H$.
Suppose that near $p \in M_{CR}$, the local CR orbits
are all of codimension 1 in $M$.
Then there exists a neighbourhood $U$ of $p$ such that
$(U \cap H) \subset \sH$, where $\sH$ is the unique Levi-flat real analytic
hypersurface in $U$ that contains $M \cap U$.
\end{lemma}

At a key point in the proof we will require the following theorem by Rossi
(See Theorem 4.7 in \cite{Rossi}).  Let $K$ be a compact subset of $X$,
$\sO(K)$ the algebra of functions holomorphic on a neighbourhood of $K$,
and $A(K)$ the closure of $\sO(K)$ in the uniform norm on $K$.
Let $M(\sA)$ be the set of peak points of the algebra $\sA$, and
$\Gamma(\sA)$ be the \v{S}ilov boundary of the algebra $\sA$.  In particular,
$p \in M(\sO(K))$ if there exists an $f \in \sO(K)$ such that
$f(p) = 1$ and $\abs{f(q)} < 1$ for $q \not= p$.

\begin{thm}[Rossi] \label{rossithm}
Suppose $K \subset X$ is compact and
$K = \bigcap_1^\infty U_k$, where $U_k$ are Stein submanifolds of $X$.
Then $M(\sO(K))$ is dense in $\Gamma(A(K))$.
\end{thm}

In the same paper (Lemma 2.9) Rossi also shows that
if $K$ is convex with respect
to the holomorphic functions on $X$, then
$K$ satisfies the
hypothesis in Theorem \ref{rossithm}, i.e. $K = \bigcap_1^\infty U_k$
for Stein submanifolds $U_k$ of $X$.

\begin{proof}[Proof of Theorem \ref{regthmtech}]
Let $X$ be a Stein manifold of dimension $N \geq 3$.
Assume that $M$ and $H$ are as in the statement
of the Theorem.  That is,
$M \subset X$ is
a compact real analytic submanifold of codimension 2,
$H$ is a connected compact $C^\infty$ Levi-flat hypersurface with boundary,
such that $\partial H = M$.

A priori, the local CR orbit of $M_{CR}$ at any particular point
can be either codimension 0, 1, or 2 in $M$.
We see from Lemma \ref{lemma:nomin}
that the local CR orbit cannot be of codimension 0 in $M$.
If some local CR orbit is of codimension 2 in $M$,
then it is
in fact a germ of a complex submanifold.  By the theorem of 
Diederich and Fornaess \cite{DF:realbnd}, no compact real analytic
submanifold contains a germ of a nontrivial complex submanifold.
If $N \geq 3$ then such a local CR orbit would in fact be nontrivial.
Hence if $M$ is a compact real analytic boundary of a Levi-flat hypersurface,
then all the local CR orbits at all points of $M_{CR}$ must be of
codimension 1 in $M$.

We will call $B \subset H$ the set of points where
$H$ is not real analytic.  First, it is obvious that $B$ is a compact
set.  By applying Lemma \ref{lemma:uniqueH}, we can see that
$H$ is real analytic, and thus extends past $M$,
near all the CR points of $M$.  Hence $B$ must be a proper closed subset of
$H$, such that $B \cap M_{CR} = \emptyset$.  Recall that $S$ is the set of CR
singular points of $M$ and $\hat{S}$ is the $\sO(X)$-convex hull of $S$.

\begin{lemma} \label{lemma:claim1}
Let $B$ and $S$ be compact subsets of $X$.
If $B \setminus \hat{S}$ is nonempty,
then there exists a point $p \in B \setminus \hat{S}$,
a neighbourhood
$U \subset X$ of $p$, and a holomorphic function $f$ defined on $U$,
such that
\begin{equation}
f(p) = 1 ~\text{ and }~ \abs{f(z)} < 1
~\text{ for all $z \in B \cap U \setminus \{p\}$}.
\end{equation}
\end{lemma}

\begin{proof}
Let us first take $\hat{B}$ be the convex hull with respect to the
holomorphic functions on $X$.  We will prove a stronger statement
than we need by showing there must exist a point
$p \in M(\sO(\hat{B})) \cap (B \setminus \hat{S})$.

The \v{S}ilov boundary $\Gamma(A(\hat{B}))$ is a subset of $B$.
Further, the $\sO(X)$-convex hull of $\Gamma(A(\hat{B}))$
is $\hat{B}$.
We assume that $B \setminus \hat{S}$ is nonempty,
therefore $\hat{B} \setminus \hat{S}$ is nonempty.
Hence, there exists a $q \in \Gamma(A(\hat{B}))$ such that
$q \in B \setminus \hat{S}$.

We apply Theorem \ref{rossithm} to
$\hat{B}$ to find that $M(\sO(\hat{B}))$ is dense in $\Gamma(A(\hat{B}))$.
Since $\hat{S}$ is closed and $q \notin \hat{S}$, there exists a point $p \in 
M(\sO(\hat{B}))$ near $q$, such that $p \notin \hat{S}$.
Since $M(\sO(\hat{B})) \subset \Gamma(A(\hat{B})) \subset B$,
we see that $p \in M(\sO(\hat{B})) \cap (B \setminus \hat{S})$.
\end{proof}

Assume for contradiction that $B \setminus \hat{S}$ was nonempty,
and that $p$ exists by Lemma \ref{lemma:claim1}.
Do note that $p \in H^o$ because $B \setminus S \subset H^o$.
Suppose that $U$ is small enough that all leafs of the Levi foliation
of $H \cap U$
are closed.
If $L \subset H \cap U$ is the closed leaf of the
Levi foliation of $H \cap U$ such that $p \in L$, then
the set $\{ \abs{f} \geq 1 \} \cap L$ has $p$ as a limit point by
the maximum principle.
Hence $L \setminus B$ also has a limit point $p$.

\begin{claim} \label{claim2}
For a perhaps smaller neighbourhood $U$ of $p$,
there exist open sets $V \subset W \subset \C^{N-1} \times \R$,
with $0 \in W$, and a
mapping $\varphi \colon W \to X$ that is
holomorphic
in the first $N-1$ variables, such that
\begin{enumerate}[(i)]
\item $\varphi(0) = p$,
\item $\varphi$ is one to one,
\item $\varphi|_V$ is real analytic and a diffeomorphism onto its image,
\item if $L$ is the leaf of the Levi foliation of $H \cap U$
through $p$, then $\varphi(V)$ contains one of the connected components of
$L\setminus B$,
\item $p \in \overline{\varphi(V)}$.
\end{enumerate}
\end{claim}

\begin{proof}
Assume that $p=0$.
Further, we can suppose that $U$ is small enough such that we can make a local
change of coordinates to assume that the leaf $L$ of the Levi foliation
through $0$ is defined by $\{ z_1 = 0 \} \cap U$.
We can assume that $U$ is perhaps even smaller, such that all
the leafs of the Levi foliation of $H \cap U$ are connected graphs
over $L$.  For any leaf $L'$ of the foliation we have
a well defined holomorphic function $\zeta \mapsto \psi_{L'}(\zeta)$, where $\zeta \in \C^{N-1}$,
whose graph over $L$ is $L'$.

Pick a point $q$ on some connected component of $L \setminus B$,
whose closure includes $p$.  Such a $q$ exists for example on the set
$\{ \abs{f} = 1 \} \cap L$.
As $H$ is real analytic near $q$,
$H$ is locally defined by $\{ \Im w_1 = 0 \}$ for some local coordinates
$w$.  We can now define $\varphi(\zeta,t) := (\zeta,\psi_{L_t}(\zeta))$,
where $L_t$ is the leaf that near $q$ is locally defined by
$\{ w_1 = t \}$.  What is left to do is to find $V$ and
show that $\varphi|_V$ is a diffeomorphism onto the right set.

Obviously, $\varphi$ is real analytic near $\varphi^{-1}(q)$.
We can cover the connected component of
$L\setminus B$ that contains $q$, by a locally finite collection
of neighbourhoods $V_j \subset L$, such that
for every point in $V_j$
a local change of coordinates as above is possible.
As we can identify $L$ with an open set in $\C^{N-1}$,
we will identify $V_j \subset L \subset \C^N$
and the corresponding set in $\C^{N-1}$.
It is not hard to see that if $V_j$ and $V_k$ overlap and
$\varphi$ is analytic on $V_j \times (-\epsilon,\epsilon)$,
then perhaps for some smaller $\epsilon'$, $\varphi$
is real analytic on $(V_j \cup V_k) \times (-\epsilon',\epsilon')$.
It is also not hard to see that if the derivative in the last variable
does not vanish on $V_j \times \{0\}$, then it does not vanish
on $V_k \times \{ 0 \}$.

Let $V = \bigcup_j V_j \times (-\epsilon_j,\epsilon_j)$.  $\varphi$ is
obviously one
to one and with a nonvanishing Jacobian on $V$, hence is a diffeomorphism onto
its image.
\end{proof}

As $\varphi|_V$ is real analytic,
we can complexify the last variable in $\varphi|_V$.  We obtain
a holomorphic mapping
$\tilde{\varphi}$ defined on an open set $\tilde{V} \subset \C^{N-1} \times \C$,
where we can think of $V \subset \tilde{V}$, such that
$\varphi|_V = \tilde{\varphi}|_V$.


Note that $0 \in \partial \tilde{V}$.
We will show that
$\tilde{V}$ is not Hartogs pseudoconvex at 0.
Let $L$ be the leaf of the Levi foliation of $H \cap U$
at $p$.  We can assume $U$ is small enough such that $L$
is closed in $H \cap U$.
Take the set
$\{ \abs{f} = 1 \} \cap L$, this is a Levi-flat real analytic
subvariety of $L$ and intersects $B$ only at $p$.
For all $\theta$ sufficiently close to 0,
the set $\{ f = e^{i\theta} \} \cap L$ is a nontrivial proper
complex subvariety of $L$ and intersects $B$ only
when $\theta = 0$ and then it intersects $B$ only at $p$.
Therefore we can find a sequence of
closed analytic discs $\Delta_j \subset L \setminus B$ such that
$\overline{\cup_j \partial \Delta_j} \subset L \setminus B$, but
$\overline{\cup_j \Delta_j} \cap B = \{ p \}$.
Hence $L \setminus B$ is not psuedoconvex at $0$ and
hence $\tilde{V}$ cannot be pseudocovex at $0$.

Since $\tilde{V}$ is not Hartogs pseudoconvex at $0$,
then $\tilde{\varphi}$ extends to a holomorphic mapping $\Phi$ defined
on a larger set, in particular in a neighbourhood of $0$.  By uniqueness
of analytic continuation along the leaves of the Levi foliation, we see that
$\Phi$ agrees with $\varphi$ on $W$ near $0$, and hence $\varphi$ is
real analytic near $0$.

Thus near $p$, $H$ is a subanalytic set.  It is also a $C^\infty$
hypersurface, and by Theorem \ref{thm:subansmooth},
$H$ is real analytic at $p$.
Hence $p \notin B$, but that would be a contradiction.
$B \setminus \hat{S}$ must have
been empty to begin with, and we are done.
\end{proof}

The technique in the proof of Theorem \ref{regthmtech} suggests other
similar results.

\begin{prop}
Suppose $H \subset X$ is a $C^\infty$ Levi-flat hypersurface (without
boundary).  Suppose $K \subset H$ is a compact set such that $H \setminus
K$ is a real analytic submanifold.  Then $H$ is a real analytic submanifold.
\end{prop}

To prove the proposition, we again let $B \subset K$ be the set where $H$
is not real analytic.  Let $S = \emptyset$
and apply Lemma \ref{lemma:claim1} to find a $p \in B$, a neighbourhood $U$ and
a peaking function $f$.
Then we follow the proof
of Theorem \ref{regthmtech} from that point on to show that $p \notin B$,
thereby showing that $B$ must have been empty to begin with.

\begin{proof}[Proof of the uniquness in Theorem \ref{regthm}]
Let $M$ be as in the statement of the theorem.
Again let $S$ be the set of CR singular points of $M$.
Suppose $H$ and $H'$ are two compact connected Levi-flat $C^\infty$
hypersurfaces with boundary $M$.  As we assume $S = \hat{S}$,
Theorem \ref{regthmtech} shows that
$H^o$ and $H'^o$ are both real analytic.
By Lemma \ref{lemma:uniqueH}, for a point $q \in M_{CR}$, there is a
neighbourhood $U_q \subset X$ of $q$
and a unique Levi-flat hypersurface $\sH$ that contains
both $H \cap U_q$ and $H' \cap U_q$.
We can assume that both $H \cap U_q$ and $H'
\cap U_q$ are connected.
There are two possibilities:
\begin{equation}
H \cap U_q = H' \cap U_q ~\text{ or }~ H \cap H' \cap U_q = M \cap U_q.
\end{equation}
In the first case,
because $H$ and $H'$ are path connected and real
analytic, unique continuation implies that $H = H'$.  Let us therefore
suppose that
$H \cap H' \cap U_q = M \cap U_q$
for all points $q \in M_{CR}$.

The compact set $H \cup H'$ is therefore such that
$(H \cup H') \setminus S$ is a Levi-flat subvariety of codimension 1
of $X \setminus S$.
Now note that $\hat{S} = S$.
If we apply Lemma \ref{lemma:claim1},
there exists a point $p \in H \cup H'$ such that
$p \notin S$, and a peaking function $f$.
That is, there exists a 
neighbourhood $V$ of $p$, and a function $f$ holomorphic in $V$
such that $\abs{f(z)} < 1$ for
$z \in V \cap (H \cup H') \setminus \{ p \}$ and $f(p) = 1$.  Simply take
$B$ to be $H \cup H'$ in Lemma \ref{lemma:claim1}.
Since
$(H \cup H') \setminus S$ is locally a union of complex hypersurfaces,
the existance of $p$ and $f$
violates the maximum principle.
\end{proof}

\section{Real algebraic extension} \label{section:algreg}

Ideally, one would want to have an extension result in the same
spirit as those in
\cite{Lebl:ext}, but for CR singular points.  The following lemma is a
small step in this direction and will be sufficient to prove
Theorem \ref{algreg}.
As the following lemma is local, note that
$H$ need not be closed.  Hence, $\partial H$ need not include all the points
in $\overline{H}$ even though $H$ is embedded.

\begin{lemma} \label{algreglocal}
Let $0\in M \subset \C^N$, $N \geq 3$,
be a real-algebraic submanifold of codimension 2,
such that $0$ is a CR singular point of $M$
and not all local CR orbits of $M_{CR}$ are codimension 2 in $M$.
Suppose there exists a connected $C^\infty$ hypersurface $H$ with boundary,
such that $H^o$ is real analytic,
$H$ is Levi-flat and $\partial H = M$.

Then there exists a local algebraic change of coordidnates near $0$
such that $M$ is locally given in $(z,w) \in \C^{N-1} \times \C$
by an equation of the form
\begin{equation}
w = \varphi(z,\bar{z}) ,
\end{equation}
where $\varphi$ is real valued.
\end{lemma}

Hence, $H$ is locally given by the equation $\Im w = 0$, and therefore
extends past $M$ near 0.
For the proof, we will require the following theorem.
\emph{Nowhere minimal} means that all local CR orbits of $M$ are 
of positive codimension in $M$.
We say that $M$ is \emph{generic} if and only if
$T_pM + J(T_pM) = T_p\C^N$ for all $p \in M$, where $J$ is the complex
structure on $\C^N$.
By $H^*$, we mean the nonsingular 
points of $H$ of top dimension (i.e. codimension 1 in $\C^N$).

\begin{thm}[See \cite{Lebl:lfnm}] \label{algcodim2thm}
Let $(M,0) \subset \C^N$ be a germ of a
real algebraic nowhere minimal generic submanifold of
codimension 2.
Then there exists an irreducible germ of a Levi-flat real algebraic 
subvariety $(H,0)$ of codimension 1, such that
for some representatives $M$ and $H$ of the germs we have
that $M \subset \overline{H^*}$.
Moreover, if not all local CR orbits of some connected representative
of $(M,0)$ are of codimension 2 in $M$, then $(H,0)$ is unique.
\end{thm}

\begin{proof}[Proof of Lemma \ref{algreglocal}]
Let $S$ be the set of CR singularities of $M$.  We note that
$\sO(\C^N)$-convex sets are precisely those that are polynomially convex.
For a proof of this fact see for example Rossi \cite{Rossi} (Lemma 2.4).
Hence we can apply Theorem \ref{regthm}
to see that $H^o$ must be real analytic.

We find a point $p \in M_{CR}$ where the CR orbit is of codimension 1.
Such a point exists by the hypothesis and by Lemma \ref{lemma:nomin}.
We then apply Theorem \ref{algcodim2thm} so that we get a real algebraic
subvariety (i.e.\@ defined by a real polynomial)
$\sH \subset \C^N$, that contains $M$.
By the uniqueness of
the germ and analytic continuation we know that $H \subset \sH$.

We can now apply Theorem \ref{thmmalvar} to know that $H$ extends
analytically past $0$, and hence extends as a real algebraic Levi-flat
hypersurface.  We know that such hypersurfaces are locally given by
an equation of the form $\Im w = 0$ in some local coordinates
$(z,w) \in \C^{N-1} \times \C$.  Since $M$ is CR singular at $0$,
we must have that $\frac{\partial}{\partial (\Re w)}$ is not
tangent to $M$ at 0, and hence
the equations
\begin{equation}
\Re w = \varphi (z, \bar{z}) ~\text{ and }~ \Im w = 0
\end{equation}
define $M$ near 0.
\end{proof}

\begin{proof}[Proof of Theorem \ref{algreg}]
Let $M$, $H$ be as in the statement of Theorem \ref{algreg}.
We first apply Theorem \ref{regthm} to ensure that $H^o$ is real analytic.
From the proof of Theorem \ref{regthmtech}
we know that all local CR orbits of $M_{CR}$
are of codimension 1 in $M$, and further that $H$ extends as a Levi-flat
hypersurface past all the CR points of $M$.  It remains to show that $H$
extends near the CR points.  We need only apply
Lemma \ref{algreglocal} as we have now satisfied all the needed hypotheses.
\end{proof}

\section{Hypersurfaces with singularities} \label{section:hypersing}

It is possible that a there might exist a singular hypersurface $H$ with
boundary $M$, even if $M$ is nonsingular.  
For example, we note that the Levi-flat subvariety $\Im (z^2 + w^2) = 0$
has an isolated singularity at 0.  If we also consider the inequality
$\abs{z}^2 + \abs{w}^2 \leq 1$, then we have a
hypersurface with a singularity at 0 and a nonsingular real analytic boundary.  This example can be generalized to $\C^N$ for $N \geq 3$ in the obvious way.

Similarly, the existence result of Dolbeault, Tomassini and
Zaitsev \cite{DTZ:CRmath} does not guarantee a nonsingular hypersurface.  An
analysis of the proof of Theorem \ref{regthmtech} allows us to formulate an
alternative statement which allows singularities.

\begin{thm} \label{regthmsing}
Let $M \subset X$
be a compact real analytic submanifold of codimension 2.
Suppose there exists a compact connected set $H$, and
a closed set $E \subset H$, such that
$H \setminus E$ is a $C^\infty$ Levi-flat hypersurface with boundary,
and $\partial (H \setminus E) = M \setminus E$.

Let $S$ be the set of CR singularities of $M$ and let
$\widehat{S \cup E}$ be the $\sO(X)$-convex hull of $S \cup E$,
then 
$H \setminus (\widehat{(S \cup E)} \cup M)$
is real analytic.
\end{thm}

In particular, if both $S$ and $E$ are a finite set,
then $H \setminus (E \cup M)$ is real analytic.

\begin{proof}
We notice that the only part of the proof of Theorem \ref{regthmtech}
that needs to
be modified is the application of Lemma \ref{lemma:claim1}.
Instead of $S$ we need to take the set $S \cup E$.
We can then apply Lemma \ref{lemma:claim1} to find a point $p \in 
B \cap \widehat{(S \cup E)}$ and a peaking function $f$.
We note that $p \in H \setminus
(M \cup \widehat{(S \cup E)})$ and derive a contradiction in the same manner as
in the proof of Theorem \ref{regthmtech}.
\end{proof}

Similarly, notice that we can modify Theorem \ref{algreg}.  Again,
we note
that the $\sO(\C^N)$-convex hull is the same as the polynomially
convex hull, and denote this convex hull of a set $K$ by
$\hat{K}$.

\begin{thm} \label{algregsing}
Let $M \subset \C^N$, $N \geq 3$,
be a compact real algebraic submanifold of codimension 2.
Suppose there exists a compact connected set $H$, and
a closed set $E \subset H$, such that
$H \setminus E$ is a $C^\infty$ Levi-flat hypersurface with boundary,
and $\partial (H \setminus E) = M \setminus E$.

Let $S$ be the set of CR singularities of $M$, if $H \setminus
\widehat{(S \cup E)}$ is connected and $M \setminus E$ is nonempty, then
there exists a real algebraic
Levi-flat subvariety $\sH \subset \C^N$ of codimension 1
such that $H \setminus \widehat{(S \cup E)} \subset \sH$.
\end{thm}

In particular, if $S$ and $E$ are finite sets then the theorem applies and
we conclude that
$H \subset \sH$.  In this case the theorem says that 
$H$ is actually a contained in a Levi flat subvariety even at points of $E$.

\begin{proof}
First we apply Theorem \ref{regthmsing} to show that
$H \setminus (\widehat{(S \cup E)} \cup M)$
is real analytic.  Since $M \setminus E$ is nonempty, then $M_{CR} \setminus
E$ is nonempty and we find a point $p \in M_{CR}$ where we can apply
Theorem \ref{algcodim2thm} to find $\sH$.
That $H \setminus \widehat{(S \cup E)} \subset \sH$ follows by analytic
continuation.
\end{proof}

\section{Examples} \label{section:examples}

To illustrate the ideas in the main theorems, we give the following examples.
First in Example \ref{example:model} we give the model example when the
theorems are true.  In Example \ref{example:uniq} we show that uniqueness
does not hold in two dimensions.  Example \ref{example:nonCR} shows that
near CR singular points, neither regularity nor uniqueness holds locally
without extra hypotheses.
Further examples of local behavior of Levi-flat hypersurfaces near a CR
boundary can be found in \cite{Lebl:ext}.

\begin{example} \label{example:model}
Suppose that $(z,w) \in \C^n \times \C$ are our coordinates and
\begin{equation}
M := \big\{ (z,w) ~\big|~ \Im w = 0, \norm{z}^2 + (\Re w)^2 = 1 \big\},
\end{equation}
where $\norm{\cdot}$ is the standard Euclidean norm on $\C^n$.  $M$
has two CR singularities, at $(0,\pm 1)$.
Then
obviously $H$ is defined by,
\begin{equation}
H := \big\{ (z,w) ~\big|~ \Im w = 0, \norm{z}^2 + (\Re w)^2 \leq 1 \big\}.
\end{equation}
It is clear that $H$ extends past its boundary as Theorem \ref{algreg}
implies.
\end{example}

\begin{example} \label{example:uniq}
The failure of uniqueness on $\C^2$ is illustrated by the following example.
Suppose that $(z,w) \in \C^2$ are our coordinates and
\begin{equation}
M := \big\{ (z,w) ~\big|~ \abs{z} = \abs{w} = 1 \big\}.
\end{equation}
There are two obvious compact Levi-flat hypersurfaces with boundary $M$,
\begin{equation}
\big\{ (z,w)  ~\big|~ \abs{z} = 1,  \abs{w} \leq 1 \big\}
\ \text{ and }\ 
\big\{ (z,w) ~\big|~ \abs{z} \leq 1,  \abs{w} = 1 \big\} .
\end{equation}
There is also a singular Levi-flat hypersurface with boundary $M$,
\begin{equation}
\big\{ (z,w) ~\big|~ \abs{z} = \abs{w},  \abs{z} \leq 1 \big\} .
\end{equation}
\end{example}

\begin{example} \label{example:nonCR}
Let us see why local regularity need not necessarily hold near CR singular
points if we do not insist on the hypersurface being compact.

First let
\begin{equation}
\varphi(x) :=
\begin{cases}
e^{-1/x} & x > 0 ,
\\
0 & x \leq 0 .
\end{cases}
\end{equation}

Let $(z,w) \in \C^2 \times \C$ be our coordinates.
Let $H$ be a Levi-flat hypersurface with boundary defined by
\begin{equation}
\begin{split}
& \Im w = \varphi(-\Re w) ,  \\
& \Re w \leq \abs{z_1}^2 + \abs{z_2}^2 .
\end{split}
\end{equation}
Note that the boundary of $H$, let us call it $M$, is defined by
\begin{equation}
\Re w = \abs{z_1}^2 + \abs{z_2}^2 \ \text{ and } \ \Im w = 0 .
\end{equation}
Outside of the origin, $M$ is a CR submanifold, where the codimension in $M$
of the CR orbits must be 1, as $M$ contains no complex analytic
subvarieties.  $M$ is in fact a real algebraic submanifold,
while $H$ clearly is not real analytic near $0$.

By multiplying $\varphi$ by a parameter we see that there are infinitely many
smooth Levi-flat hypersurfaces with boundary $M$.  On the other hand, because
the CR orbits are of codimension 1, Lemma \ref{lemma:uniqueH} and analytic
continuation shows that there are only two connected real analytic Levi-flat
hypersurfaces with boundary $M$.
One is defined by
$\Im w = 0$ and $\Re w \leq \abs{z_1}^2 + \abs{z_2}^2$,
and the other by
$\Im w = 0$ and $\Re w \geq \abs{z_1}^2 + \abs{z_2}^2$.
\end{example}


\def\MR#1{\relax\ifhmode\unskip\spacefactor3000 \space\fi%
  \href{http://www.ams.org/mathscinet-getitem?mr=#1}{MR#1}}

\end{document}